\documentclass{amsart}
\usepackage{amsmath, amssymb, latexsym}
\title{Face Vectors of Flag Complexes}
\author{Andy Frohmader}
\address{Department of Mathematics, University of
  Washington, Seattle, WA 98195-4350}
\email{frohmade@math.washington.edu}

\newtheorem{theorem}{Theorem}[section]

\newtheorem{lemma}[theorem]{Lemma}
\newtheorem{definition}[theorem]{Definition}

\newcommand{\link}[2]{\ensuremath{\textrm{lk}_{#1}({#2})}}

\begin{document}

\begin{abstract}
A conjecture of Kalai and Eckhoff that the face vector of an
arbitrary flag complex is also the face vector of some particular
balanced complex is verified.
\end{abstract}

\date{May 8}

\maketitle

\section {Introduction}

We begin by introducing the main result.  Precise definitions and
statements of some related theorems are deferred to later
sections.

The main object of our study is the class of flag complexes.  A
simplicial complex is a \textit{flag complex} if all of its
minimal non-faces are two element sets. Equivalently, if all of
the edges of a potential face of a flag complex are in the
complex, then that face must also be in the complex.

Flag complexes are closely related to graphs.  Given a graph $G$,
define its \textit{clique complex} $C = C(G)$ as the simplicial
complex whose vertex set is the vertex set of $G$, and whose faces
are the cliques of $G$.  The clique complex of any graph is itself
a flag complex, as for a subset of vertices of a graph to not form
a clique, two of them must not form an edge. Conversely, any flag
complex is the clique complex of its 1-skeleton.

The Kruskal-Katona theorem \cite{kruskal, katona} classifies the
face vectors of simplicial complexes as being precisely the
integer vectors whose coordinates satisfy some particular bounds.
The graphs of the ``rev-lex'' complexes which attain these bounds
invariably have a clique on all but one of the vertices of the
complex, and sometimes even on all of the vertices.

Since the bounds of the Kruskal-Katona theorem hold for all
simplicial complexes, they must in particular hold for flag
complexes.  We might expect that flag complexes which do not have
a face on most of the vertices of the complex will not come that
close to attaining the bounds of the Kruskal-Katona theorem.

One way to force tighter bounds on face numbers is by requiring
the graph of the complex to have a chromatic number much smaller
than the number of vertices.  The face vectors of simplicial
complexes of a given chromatic number were classified by Frankl,
F\"{u}redi, and Kalai \cite{balanced}.

Kalai (unpublished; see \cite[p.~100]{greenbook}) and Eckhoff
\cite{mainconj} independently conjectured that if the largest face
of a flag complex contains $r$ vertices, then it must satisfy the
known bounds (see \cite{balanced}) for complexes of chromatic
number $r$, even though the flag complex may have chromatic number
much larger than $r$. We prove their conjecture.

\begin{theorem}
For any flag complex $C$, there is a balanced complex $C'$ with
the same face vector as $C$.\label{maintheorem}
\end{theorem}

Our proof is constructive.  The Frankl-F\"{u}redi-Kalai
\cite{balanced} theorem states that an integer vector is the face
vector of a balanced complex if and only if it is the face vector
of a colored ``rev-lex'' complex.  This happens if and only if it
satisfies certain bounds on consecutive face numbers.  Given a
flag complex, for each $i$, we construct a colored ``rev-lex''
complex with the same number of $i$-faces and $(i+1)$-faces as the
flag complex, thus showing that all the bounds are satisfied.

The structure of the paper is as follows.  Section 2 contains
basic facts and definitions related to simplicial complexes.  In
Section 3, we discuss the Kruskal-Katona theorem and the
Frankl-F\"{u}redi-Kalai theorem, and lay the foundation for our
proof.  Finally, Section 4 gives our proof of the Kalai-Eckhoff
conjecture.

\section {Preliminaries on simplicial complexes}

In this section, we discuss some basic definitions related to
simplicial complexes.

Recall that a \textit{simplicial complex} $\Delta$ on a vertex set
$V$ is a collection of subsets of $V$ such that, (i) for every $v
\in V$, $\{v\} \in \Delta$ and (ii) for every $B \in \Delta$, if
$A \subset B$, then $A \in \Delta$.  The elements of $\Delta$ are
called \textit{faces}.  The maximal faces (under inclusion) are
called \textit{facets}.

For a face $F$ of a simpicial complex $\Delta$, the
\textit{dimension} of $F$ is defined as dim $F = |F|-1$.  The
dimension of $\Delta$, dim $\Delta$, is defined as the maximum
dimension of the faces of $\Delta$.  A complex $\Delta$ is
\textit{pure} if all of its facets are of the same dimension.

The $i$-skeleton of a simplicial complex $\Delta$ is the
collection of all faces of $\Delta$ of dimension $\leq i$.  In
particular, the 1-skeleton of $\Delta$ is its underlying graph.

It is sometimes useful in inductive proofs to consider certain
subcomplexes of a given simplicial complex, such as its links.

\begin{definition}
\textup{Let $\Delta$ be a simplicial complex and $F \in \Delta$.
 The \textit{link} of $F$, \link{\Delta}{F}, is defined as}
$$\link{\Delta}{F} := \{G \in \Delta \ |\ F \cap G = \emptyset, F \cup G \in \Delta\}.$$
\end{definition}

The link of a face of a simplicial complex is itself a simplicial
complex.  It will be convenient to define the notion of a link of
a vertex of a graph.

\begin{definition}
The link of a vertex $v$ in a graph $G$, denoted \link{G}{v}, is
the induced subgraph of $G$ on all vertices adjacent to $v$.
\end{definition}

Note that \link{G}{v} coincides with the 1-skeleton of the link of
$v$ in the clique complex of $G$.

Next we discuss a special class of simplical complexes known as
flag complexes.

\begin{definition}
\textup{A simplicial complex $\Delta$ on a vertex set $V$ is a
\textit{flag complex} if all of its minimal non-faces are two
element sets.  A non-face of $\Delta$ is a subset $A \subseteq V$
such that $A \not \in \Delta$.  A non-face $A$ is minimal if, for
all proper subsets $B \subset A$, $B \in \Delta$.}
\end{definition}

In the following, we refer to the chromatic number of a simplicial
complex as the chromatic number of its 1-skeleton in the usual
graph theoretic sense.

We also need the notion of a balanced complex, as introduced and
studied in \cite{gencolor}.

\begin{definition}
\textup{A simplicial complex $\Delta$ of dimension $d-1$ is
\textit{balanced} if it has chromatic number $d$.}
\end{definition}

Note that the chromatic number of a simpicial complex of dimension
$d-1$ must be at least $d$, as it has some face with $d$ vertices,
all of which are adjacent, so coloring that face takes $d$ colors.
A balanced complex is then one whose chromatic number is no larger
than it has to be.

Not all simplicial complexes are balanced complexes.  For example,
a pentagon (five vertices, five edges, and one empty face) is not
a balanced complex, because it has chromatic number three but
dimension only one.

In this paper, we study the face numbers of flag complexes.

\begin{definition}
\textup{The \textit{$i$-th face number} of a simplicial complex
$C$, denoted $c_i(C)$ is the number of faces in $C$ containing $i$
vertices. These are also called \textit{\mbox{$i$-faces}} of $C$.
If dim $C = d-1$, the \textit{face vector} of $C$ is the vector
$$c(C) = (c_0(C), c_1(C), \dots , c_d(C)).$$}
\end{definition}

In particular, for any non-empty complex $C$, we have $c_0(C) =
1$, as there is a unique empty set of vertices, and it is a face
of $C$.

Since flag complexes are the same as clique complexes of graphs,
it is sometimes convenient to talk about face numbers in the
language of graphs.

\begin{definition}
\textup{The \textit{$i$-th face number} of a graph is the $i$-th
face number of its clique complex.  Likewise, the \textit{clique
vector} of a graph is the face vector of its clique complex.}
\end{definition}

The face numbers defined here are shifted by one from what is
often used for simplicial complexes.  This is done because we are
primarily concerned with flag complexes, or equivalently, clique
complexes of graphs, where it is more natural to index $i$ as the
number of vertices in a clique of the graph, following Eckhoff
\cite{eckhoff}.

The graph concept corresponding to the dimension of a simplicial
complex is the clique number.

\begin{definition}
\textup{The \textit{clique number} of a graph is the number of
vertices in its largest clique.}
\end{definition}

Note that the clique number of a graph is one larger than the
dimension of its clique complex.

\section {The Kruskal-Katona and Frankl-F\"{u}redi-Kalai theorems}

For the general case of simplicial complexes, the question of
which face vectors are possible is answered by the Kruskal-Katona
theorem \cite{kruskal, katona}.  Stating the theorem requires the
following lemma.

\begin{lemma}
Given any positive integers $m$ and $k$, there is a unique $s$ and
unique $n_k > n_{k-1} > \dots > n_{k-s} \geq k-s > 0$ such that
$$m = {n_k \choose k} + {n_{k-1}\choose k-1} + \dots + {n_{k-s}
\choose k-s}.$$
\end{lemma}

The representation described in the lemma is called the
$k$-canonical representation of $m$.

\begin{theorem}
[Kruskal-Katona]  For a simplicial complex $C$, let $$m = c_k(C) =
{n_k \choose k} + {n_{k-1}\choose k-1} + \dots + {n_{k-s} \choose
k-s}$$ be the $k$-canonical representation of $m$. Then
$$c_{k+1}(C) \leq {n_k \choose k+1} + {n_{k-1}\choose k} + \dots +
{n_{k-s} \choose k-s+1}.$$ Furthermore, given a vector $(1, c_1,
c_2, \dots, c_t)$ which satisfies this bound for all $1 \leq k <
t$, there is some complex that has this vector as its face
vector.\label{kruskat}
\end{theorem}

To construct the complexes which demonstrate that the bound of the
Kruskal-Katona theorem is attained, we need the
reverse-lexicographic (``rev-lex") order.  To define the rev-lex
order of $i$-faces of a simplicial complex on $n$ vertices, we
start by labelling the vertices $1, 2, \dots$.  Let $\mathbb{N}$
be the natural numbers, let $A$ and $B$ be distinct subsets of
$\mathbb{N}$ with $|A| = |B| = i$, and let $A \nabla B$ be the
symmetric difference of $A$ and $B$.

\begin{definition}
\textup{For $A, B \subset \mathbb{N}$ with $|A| = |B|$, we say
that $A$ precedes $B$ in the rev-lex order if max$(A \nabla B) \in
B$, and $B$ precedes $A$ otherwise.}
\end{definition}

For example, $\{2, 3, 5\}$ precedes $\{1, 4, 5\}$, as 3 is less
than 4, but $\{3, 4, 5\}$ precedes $\{1, 2, 6\}$.

\begin{definition}
\textup{The \textit{rev-lex complex on $m$ $i$-faces} is the pure
complex whose facets are the first $m$ $i$-sets possible in
rev-lex order.  This complex is denoted $C_i(m)$.}
\end{definition}

We can also specify more than one number in the face vector.  For
two sequences $i_1 < \dots < i_r$ and $(m_1, \dots, m_r)$, let $$C
= C_{i_1}(m_1) \cup C_{i_2}(m_2) \cup \dots \cup C_{i_r}(m_r).$$
The standard way to prove the Kruskal-Katona theorem involves
showing that if the numbers $m_1, \dots, m_r$ satisfy the bounds
of the theorem, then the complex $C$ has exactly $m_j$ $i_j$-faces
for all $j \leq r$ and no more.  In this case, we refer to $C$ as
the rev-lex complex on $m_1$ $i_1$-faces, \dots, $m_r$
$i_r$-faces.

For example, if the complex $C$ has ${9 \choose 3} + {6 \choose 2}
= 99$ 3-faces, then the Kruskal-Katona theorem says that it can
have at most ${9 \choose 4} + {6 \choose 3} = 146$ 4-faces. The
rev-lex complex on 99 3-faces and 146 4-faces gives an example
showing that this bound is attained.

The 1-skeleton of the rev-lex complex that gives the example for
the existence part of the Kruskal-Katona theorem always has as
large of a clique as is possible without exceeding the number of
edges allowed, as well as a chromatic number of either the number
of non-isolated vertices or one less than this.  It turns out that
if we require a much smaller chromatic number, we can get a much
smaller bound.  To take an extreme example, if $c_3(C) = 1140$,
then the Kruskal-Katona theorem requires that $c_4(C) \leq 4845$.
But if we require the complex $C$ to be 3-colorable, then we
trivially cannot have any faces on four vertices, and $c_4(C) =
0$.

We could ask what face vectors occur for $r$-colorable complexes
for a given $r$.  This was solved by Frankl, F\"{u}redi, and Kalai
\cite{balanced}.  In order to explain their result, we need the
concept of a Tur\'{a}n graph.

\begin{definition}
\textup{The Tur\'{a}n graph $T_{n,r}$ is the graph obtained by
partitioning $n$ vertices into $r$ parts as evenly as possible,
and making two vertices adjacent exactly if they are not in the
same part.  Define ${n \choose k}_r$ to be the number of
$k$-cliques of the graph $T_{n,r}$.}
\end{definition}

The structure of the Frankl-F\"{u}redi-Kalai theorem
\cite{balanced} is similar to that of the Kruskal-Katona theorem,
beginning with a canonical representation of the number of faces.

\begin{lemma}
Given positive integers m, k, and r with $r\geq k$, there are
unique $s$, $n_k$, $n_{k-1}$, \dots , $n_{k-s}$ such that $$m =
{n_k\choose k}_r + {n_{k-1}\choose k-1}_{r-1} + \dots +
{n_{k-s}\choose k-s}_{r-s},$$ $n_{k-i}-\big\lfloor{n_{k-i}\over
r-i}\big\rfloor > n_{k-i-1}$ for all $0\leq i < s,$ and
$n_{k-s}\geq k-s > 0$.
\end{lemma}

This expression is called the $(k, r)$-canonical representation of
$m$.

\begin{theorem}
[Frankl-F\"{u}redi-Kalai] For an $r$-colorable complex $C$, let
$$m = c_k(C) = {n_k\choose k}_r + {n_{k-1}\choose k-1}_{r-1} +
\dots + {n_{k-s}\choose k-s}_{r-s}$$ be the $(k, r)$-canonical
representation of $m$. Then $$c_{k+1}(C) \leq {n_k\choose k+1}_r +
{n_{k-1}\choose k}_{r-1} + \dots + {n_{k-s}\choose k-s+1}_{r-s}.$$
Furthermore, given a vector $(1, c_1, c_2, \dots c_t)$ which
satisfies this bound for all $1 \leq k < t$, there is some
$r$-colorable complex that has this vector as its face vector.
\label{coloredKK}
\end{theorem}

The examples which show that this bound is sharp come from a
colored equivalent of the rev-lex complexes of the Kruskal-Katona
theorem.

\begin{definition}
\textup{A subset $A \subset \mathbb{N}$ is
\textit{$r$-permissible} if, for any two $a, b \in A$, $r$ does
not divide $a-b$.  The \textit{$r$-colored rev-lex complex on $m$
$i$-faces} is the pure complex whose facets are the first $m$
$r$-permissible $i$-sets in rev-lex order.  This complex is
denoted $C_i^r(m)$.}
\end{definition}

The complex $C_i^r(n)$ is $r$-colorable because we can color all
vertices which are $i$ modulo $r$ with color $i$.

As with the uncolored case, we can define a rev-lex complex with
specified face numbers of more than one dimension.  For two
sequences $i_1 < \dots < i_s$ and $(m_1, \dots, m_s)$, let $C =
C_{i_1}^r(m_1) \cup C_{i_2}^r(m_2) \cup \dots \cup
C_{i_s}^r(m_s)$.  The proof of Theorem~\ref{coloredKK} involves
showing that if the numbers $m_1, \dots, m_r$ satisfy the bounds
of the theorem, then the complex $C$ has exactly $m_j$ $i_j$-faces
and no more. In this case, we refer to $C$ as the $r$-colored
rev-lex complex on $m_1$ $i_1$-faces, \dots, $m_r$ $i_r$-faces.
This complex is likewise $r$-colorable with one color for each
value modulo $r$.

In the case of flag complexes, the face numbers of the complex
must still follow the bounds imposed by the chromatic number by
Theorem~\ref{coloredKK}. Still, there are graphs whose clique
number is far smaller than the chromatic number, and having no
large cliques seems to force tighter restrictions on the clique
vector than the chromatic number alone. In particular, given a
graph $G$ of clique number $n$, we must have $c_i(G) = 0$ for all
$i > n$, while the bound from the chromatic number and
Theorem~\ref{coloredKK} may be rather large. Note that the
chromatic number must be at least the size of the largest clique,
as any two vertices in a maximum size clique must have different
colors.

It has been conjectured by Kalai (unpublished) and Eckhoff
\cite{mainconj} that, given a graph $G$ with clique number $r$,
there is an $r$-colorable complex with exactly the same face
numbers as the clique complex of the graph.  Their conjecture
generalizes the classical Tur\'{a}n theorem from graph theory,
which states that among all triangle-free graphs on $n$ vertices,
the Tur\'{a}n graph $T_{n,2}$ has the most edges \cite{turan}. The
goal of the following section is to verify
Theorem~\ref{maintheorem}, proving their conjecture.

\section {Proof of the Kalai-Eckhoff conjecture}

Fix a graph $G$ with $c_{r+1}(G) = 0$ and fix $k \geq 0$.  We
start by showing that there is an $r$-colorable complex $C$ with
$c_k(G) = c_k(C)$ and $c_{k+1}(G) = c_{k+1}(C)$ (see
Lemma~\ref{mainlemma} below).

The case $k = 1$ of the lemma is given by Tur\'{a}n's theorem
\cite{turan}.  It was generalized by Zykov \cite{zykov} to state
that if $G$ is a graph on $n$ vertices of chromatic number $r$,
then $c_i(G) \leq {n \choose i}_r$.  The case $k = 2$ was proven
by Eckhoff \cite{eckhoff0}.  A subsequent paper of Eckhoff
\cite{eckhoff} established a bound on $c_i(G)$ in terms of
$c_2(G)$ for all $2 \leq i$.  All of these results are special
cases of our Theorem~\ref{restate} and proven independently below.

\begin{lemma}
If $G$ is a graph with $c_{r+1}(G) = 0$ and $k$ is a nonnegative
integer, then there is some $r$-colorable complex $C$ with $c_k(C)
= c_k(G)$ and $c_{k+1}(C) = c_{k+1}(G)$. \label{mainlemma}
\end{lemma}

Proof.  We use induction on $k$.  For the base case, if $k = 0$,
take $C$ to be a complex with the same number of vertices as $G$,
no edges, and all vertices the same color.

Otherwise, assume that the lemma holds for $k-1$, and we need to
prove it for $k$.  The approach for this is to use induction on
$c_{k+1}(G)$.  For the base case, if $c_{k+1}(G) = 0$, then there
is trivially some $r$-colorable complex $C$ with $c_k(C) = c_k(G)$
and $c_{k+1}(C) = 0$.

For the inductive step, suppose that $c_{k+1}(G) > 0$. Let $v_0$
be the vertex of $G$ contained in the most cliques of $k+1$
vertices; in case of a tie, arbitrarily pick some vertex tied for
the most to label $v_0$.  Let the vertices of $G$ not adjacent to
$v_0$ be $v_1, v_2, \dots v_s$.

Given a graph $G$ and a vertex $v$, there is a bijection between
$k$-cliques of \link{G}{v} and $(k+1)$-cliques of $G$ containing
$v$, where a $k$-clique of \link{G}{v} corresponds to the
$(k+1)$-clique of $G$ containing the $k$ vertices of the
$k$-clique of \link{G}{v} together with $v$.  Then the number of
$(k+1)$-cliques of $G$ containing $v$ is $c_k(\link{G}{v})$.  In
particular, the choice of $v_0$ gives $c_k(\link{G}{v_0}) \geq
c_k(\link{G}{v'})$ for every vertex $v' \in G$.

Define graphs $G_0, G_1, \dots, G_{s+1}$ by setting $G_{i+1} = G -
\{v_0, v_1, \dots v_i\}$ for $0 \leq i \leq s$ and $G_0 = G$.
Clearly, $G = G_0 \supset G_1 \supset \dots \supset G_{s+1}$.
Further, $G_{s+1}$ is the induced subgraph on the vertices
adjacent to $v_0$, which is \link{G}{v_0}.

Since $c_{r+1}(G) = 0$, $c_r(\link{G}{v_0}) = 0$, for otherwise,
the $r$ vertices of an $r$-clique of \link{G}{v_0} together with
$v_0$ would form an $(r+1)$-clique of $G$. Then $c_r(G_{s+1}) =
0$. Further, since $c_{k+1}(G) > 0$, and $v_0$ is contained in the
most $(k+1)$-cliques of any vertex of $G$, $v_0$ is contained in
at least one $(k+1)$-clique, and so $c_k(\link{G}{v_0}) > 0$.
Since $v$ is contained in at least one $(k+1)$-clique of $G$, we
have $c_{k+1}(G_{s+1}) < c_{k+1}(G)$.

Then by the second inductive hypothesis, there is some
$(r-1)$-colorable complex $C_{s+1}$ such that $c_k(C_{s+1}) =
c_k(G_{s+1})$ and $c_{k+1}(C_{s+1}) = c_{k+1}(G_{s+1})$. Since
given any $(r-1)$-colorable complex, there is an $(r-1)$-colorable
rev-lex complex with the same face numbers, we can take $C_{s+1}$
to be a rev-lex complex. Further, since $c_{k+1}(C_{s+1})$ and
$c_k(C_{s+1})$ only force a lower bound on $c_{k-1}(C_{s+1})$, but
not an upper bound, we can take $c_{k-1}(C_{s+1}) \geq
c_{k-1}(G)$.

Let $c_k(\link{G_i}{v_i}) = a_i$ and $c_{k-1}(\link{G_i}{v_i}) =
b_i$. Since $G_{i+1} = G_i - v_i$, $c_{k+1}(G_i) -
c_{k+1}(G_{i+1}) = a_i$ and $c_k(G_i) - c_k(G_{i+1}) = b_i$.  We
have $c_k(\link{G}{v_0}) \geq c_k(\link{G}{v_i})$ by the choice of
$v_0$. We also have $c_k(\link{G}{v_i}) \geq c_k(\link{G_i}{v_i})$
since $G_i \subset G$. Thus $$c_k(C_{s+1}) = c_k(G_{s+1}) =
c_k(\link{G}{v_0}) \geq c_k(\link{G}{v_i}) \geq
c_k(\link{G_i}{v_i}) = a_i.$$

Given an $r$-colored complex $C_{i+1}$ such that $c_{k+1}(C_{i+1})
= c_{k+1}(G_{i+1})$, $c_k(C_{i+1}) = c_k(G_{i+1})$, and the
induced subcomplex of $C_{i+1}$ on the vertices of the first $r-1$
colors is isomorphic to $C_{s+1}$, we want to construct a complex
$C_i$ such that $c_{k+1}(C_i) = c_{k+1}(G_i)$, $c_k(C_i) =
c_k(G_i)$, and the induced subcomplex of $C_i$ on the vertices of
the first $r-1$ colors is isomorphic to $C_{s+1}$.

Construct $C_i$ from $C_{i+1}$ by adding a new vertex $v'_i$ of
color $r$.  Let the $(k+1)$-faces containing $v'_i$ consist of
each of the first $a_i$ $k$-faces in rev-lex order of $C_{s+1}$
together with $v'_i$, and let the $k$-faces containing $v'_i$
consist of each of the first $b_i$ $(k-1)$-faces in rev-lex order
of $C_{s+1}$ together with $v'_i$.

If this construction can be done, then $c_{k+1}(C_i)$ is the
number of $(k+1)$-faces of $C_i$ containing $v'_i$ plus the number
of $(k+1)$-faces of $C_i$ not containing $v'_i$, which are $a_i$
and $c_{k+1}(C_{i+1})$, respectively.  Then
$$c_{k+1}(C_i) = c_{k+1}(C_{i+1}) + a_i = c_{k+1}(G_{i+1}) + a_i =
c_{k+1}(G_i).$$ Likewise, we have $$c_k(C_i) = c_k(C_{i+1}) + b_i
= c_k(G_{i+1}) + b_i = c_k(G_i).$$  Further, it is clear from the
construction that the induced subcomplex on vertices of the first
$r-1$ colors is unchanged from $C_{i+1}$, and hence is isomorphic
to $C_{s+1}$.

In order to show that the construction is possible, we must show
that $c_k(C_{s+1}) \geq a_i$ and $c_{k-1}(C_{s+1}) \geq b_i$, and
that it is possible for an $(r-1)$-colored complex $C$ to have
exactly $c_k(C) = a_i$ and $c_{k-1}(C) = b_i$.  For the first of
these, we have already shown that $c_k(C_{s+1}) \geq a_i$.

For the second, $c_{k-1}(\link{G_i}{v_i}) = b_i$.  But
$\link{G_i}{v_i} \subset G_i \subset G$, so $$b_i =
c_{k-1}(\link{G_i}{v_i}) \leq c_{k-1}(G_i) \leq c_{k-1}(G) \leq
c_{k-1}(C_{s+1}).$$

For the third, since $G_i \subset G$, we have $c_{r+1}(G_i) \leq
c_{r+1}(G) = 0$, and so $c_{r+1}(G_i) = 0$.  Then
$c_r(\link{G_i}{v_i}) = 0$.  We also have $c_k(\link{G_i}{v_i}) =
a_i$ and $c_{k-1}(\link{G_i}{v_i}) = b_i$ by the definitions of
$a_i$ and $b_i$.  Then by the first inductive hypothesis, there is
some $(r-1)$-colored complex $C'_i$ such that $c_k(C'_i) = a_i$
and $c_{k-1}(C'_i) = b_i$.  Then we can take $C'_i$ to be the
$(r-1)$-colored rev-lex complex with $c_k(C'_i) = a_i$ and
$c_{k-1}(C'_i) = b_i$. Since $C_{s+1}$ is an $(r-1)$-colored
rev-lex complex with $c_k(C_{s+1}) \geq a_i$ and $c_{k-1}(C_{s+1})
\geq b_i$, $C'_i \subset C_{s+1}$, and we can choose the link of
$v'_i$ in $C_i$ to be $C'_i$.

We can repeat this construction for each $0 \leq i \leq s$ to
start with $C_{s+1}$, then construct $C_s$, then $C_{s-1}$, and so
forth, until we have an $r$-colored complex $C_0$ such that
$c_k(C_0) = c_k(G)$ and $c_{k+1}(C_0) = c_{k+1}(G)$.  This
completes the inductive step for the induction on $c_{k+1}(G)$,
which in turn completes the inductive step for the induction on
$k$. $\Box$

We are now ready to prove the result which immediately implies
Theorem~\ref{maintheorem}, and hence establish the Kalai-Eckhoff
conjecture, by taking $r$ to be the clique number of $G$.

\begin{theorem}
For every graph $G$ with $c_{r+1}(G) = 0$, there is an
$r$-colorable complex $C$ such that $c_i(C) = c_i(G)$ for all
$i$.\label{restate}
\end{theorem}

Proof.  By Lemma~\ref{mainlemma}, we can pick an $r$-colored
complex $C_i$ such that $c_i(C_i) = c_i(G)$ and $c_{i+1}(C_i) =
c_{i+1}(G)$ for all $i \geq 1$.  By Theorem~\ref{coloredKK}, we
can take $C_i$ to be the rev-lex complex on $c_i(G)$ $i$-faces and
$c_{i+1}(G)$ $(i+1)$-faces, and then $\cup_{i=1}^{r} C_i$ will
have the desired face numbers. $\Box$

\textit{Acknowledgements.} I would like to thank my thesis advisor
Isabella Novik for suggesting this problem, and for her many
useful discussions on solving the problem and writing an article.

\end{document}